\documentclass[letterpaper]{article} 
\usepackage{aaai2026}  
\usepackage{times}  
\usepackage{helvet}  
\usepackage{courier}  
\usepackage[hyphens]{url}  
\usepackage{graphicx} 
\usepackage{amsmath}

\usepackage{natbib}  
\usepackage{caption} 
\usepackage{algorithm}
\usepackage{algorithmic}
\usepackage{newfloat}
\usepackage{listings}
\DeclareCaptionStyle{ruled}{labelfont=normalfont,labelsep=colon,strut=off} 
\floatstyle{ruled}
\newfloat{listing}{tb}{lst}{}
\floatname{listing}{Listing}
\title{A Comparison of Reinforcement Learning and \\ Optimal Control Methods for Path Planning}
\author{
    Qiang Le\textsuperscript{\rm 1},%
    Yaguang Yang\textsuperscript{\rm 1},
    Isaac E. Weintraub\textsuperscript{\rm 2}
}
\affiliations{
    \textsuperscript{\rm 1}Department of Electrical and Computer Engineering, Hampton University\\
    Hampton, USA\\
    qiang.le@hamptonu.edu, yaguang.yang@hamptonu.edu \\
   \textsuperscript{\rm 2}Air Warfare Directorate, Air Force Research Laboratory\\
Wright-Patterson, USA\\
isaac.weintraub.1@afrl.af.mil
}
\usepackage{bibentry}

\begin{document}

\maketitle

\begingroup
\renewcommand{\thefootnote}{}%
\endgroup

\begin{abstract}
Path-planning for autonomous vehicles in threat-laden environments is a fundamental challenge. While traditional optimal control methods can find ideal paths, the computational time is often too slow for real-time decision-making. To solve this challenge, we propose a method based on Deep Deterministic Policy Gradient (DDPG) and model the threat as a simple, circular `no-go' zone. A mission failure is claimed if the vehicle enters this `no-go' zone at any time or does not reach a neighborhood of the destination. The DDPG agent is trained to learn a direct mapping from its current state (position and velocity) to a series of feasible actions that guide the agent to safely reach its goal. A reward function and two neural networks, critic and actor, are used to describe the environment and guide the control efforts. The DDPG trains the agent to find the largest possible set of starting points (``feasible set'') wherein a safe path to the goal is guaranteed. This provides critical information for mission planning, showing beforehand whether a task is achievable from a given starting point, assisting pre-mission planning activities. The approach is validated in simulation. A comparison between the DDPG method and a traditional optimal control (pseudo-spectral) method is carried out. The results show that the learning-based agent may produce effective paths while being significantly faster, making it a better fit for real-time applications. However, there are areas (``infeasible set'') where the DDPG agent cannot find paths to the destination, and the paths in the feasible set may not be optimal. These preliminary results guide our future research: (1) improve the reward function to enlarge the DDPG feasible set, (2) examine the feasible set obtained by the pseudo-spectral method, and (3) investigate the arc-search IPM method for the path planning problem.
\end{abstract}

\begin{links}
    \link{Code}{Available-from-the-authors}
    \link{Datasets}{Available-from-the-authors}
    \link{Extended version}{Available-from-the-authors}
\end{links}

\section{Introduction}
Path planning in the presence of obstacles has many applications in various scenarios, for example, in automatic driving \cite{hctch18}, space debris avoidance \cite{msml11}, avoidance of engagement zones \cite{dzwv23,vw24,wvchf22}, robotics \cite{gc13}, video games \cite{ask15}, geographical information systems \cite{bg07}, printed circuit board routing \cite{zwlw22}, and practical VLSI systems \cite{up23}, among others. Therefore, a large body of research has focused on this problem; readers may find a lot of information from the survey papers \cite{gc13,ksds21,wwn20,zm18} and references therein. Several techniques, such as optimal control (pseudo-spectral and model predictive control methods) \cite{herber15}, machine learning \cite{hwmz20}, and Voronoi diagrams \cite{nstj19}, have been used to solve the problem. However, little work has been done to compare the advantages and disadvantages of the different methods.

In this paper, we investigated a DDPG-based machine learning method for path planning in the presence of obstacles. Our main purpose is to have a preliminary idea about the advantages and drawbacks of the two different methodologies: the traditional optimal control/optimization-based method and the machine learning-based method, and we hope this comparison will guide our future research direction. For this reason, the DDPG-based method is implemented in Matlab. The preliminary results are compared to those obtained in a Matlab implementation of the pseudo-spectral method \cite{lw26}. This work led us to believe that while the DDPG-based machine learning method takes significantly more computational time (several minutes) in training, a trained agent normally finds a feasible path much faster (at least an order of magnitude) than the pseudo-spectral method, which is important for real-time applications. On the other hand, the pseudo-spectral method can normally find an optimal solution while the DDPG-based machine learning method normally cannot. By using a ``feasibility map,'' we noticed that the DDPG-based method in its current form cannot find a feasible path for every starting point. Our extensive experience made us believe that this can be improved if we refine the reward function and adapt a smart initial heading. This will be the next step, and comprehensive research will follow.

\section{Related Work}
A foundational reinforcement learning algorithm is Q-learning, in which the agent updates the Q-table that stores the expected discounted long-term accumulative reward for a given observation and action. Q-learning is able to find the optimal path given an initial position, the final position, and the static no-go zones while the number of actions is small, for example, the 8 directional actions: E, W, N, S, SW, SE, NE, and NW. Because Q-learning is able to find the optimal path and is easy to use, it has been considered for solving a relatively simple path planning problem in the presence of obstacles \cite{cx18,mh20,prpplf24,syc23,wyl22}.

As the grid-based environment becomes large, the size of the Q-table explodes and becomes impossible to update. To overcome this difficulty, an idea that couples Q-learning with deep learning or a neural network is proposed. The new framework, known as DQN, is able to extract the features of the environment, and therefore handle the continuous environment. Because of this improvement, DQN has been used for the path planning problem where the state space is continuous \cite{gzlsx23,hqll23,nkm23,xyj25,ylp20}.

Although the DQN method allows infinitely many (or continuous) states, the action space for the DQN agent is still discrete. To have a method that can handle both continuous states and continuous actions, DDPG was developed by a group of Google engineers \cite{slhdwr14,lillicrap2015continuous} based on Q-learning and Deep Q-network (DQN). It uses two neural networks: one for continuous Q-value estimation and one for continuous action. Therefore, this reinforcement learning method is highly scalable. In addition, DDPG has the ability to learn deterministic policies. As a result, many recent path planning works follow this direction, for example, \cite{asz24,gylwc23,rhv23,whlm23,wges22,xcayz22}.

\section{Deep Reinforcement Learning Approach}\label{DDPGmethod}

\subsection{Actor and critic neural networks}
We select the DDPG approach to tackle the path-finding control problem given its suitability for continuous state and action spaces. Building on concepts from Q-learning, DDPG combines an actor network and a critic network to explore continuous action spaces\cite{lillicrap2015continuous}. The output of the actor network is the agent action $A$ using the existing network parameters $\theta_A$ for the given state $S$:
\begin{equation}
A=\pi(S;\theta_A),
\label{policy}
\end{equation}
where $\pi$ is the policy. The output of the critic network is the Q value for the given state $S$ and the action $A$ based on the critic network parameters $\theta_C$:
\begin{equation}
Q(S,A; \theta_C).
\label{cNetwork}
\end{equation}

\begin{figure}[htbp]
	\centering
	\includegraphics[width=0.48\textwidth]{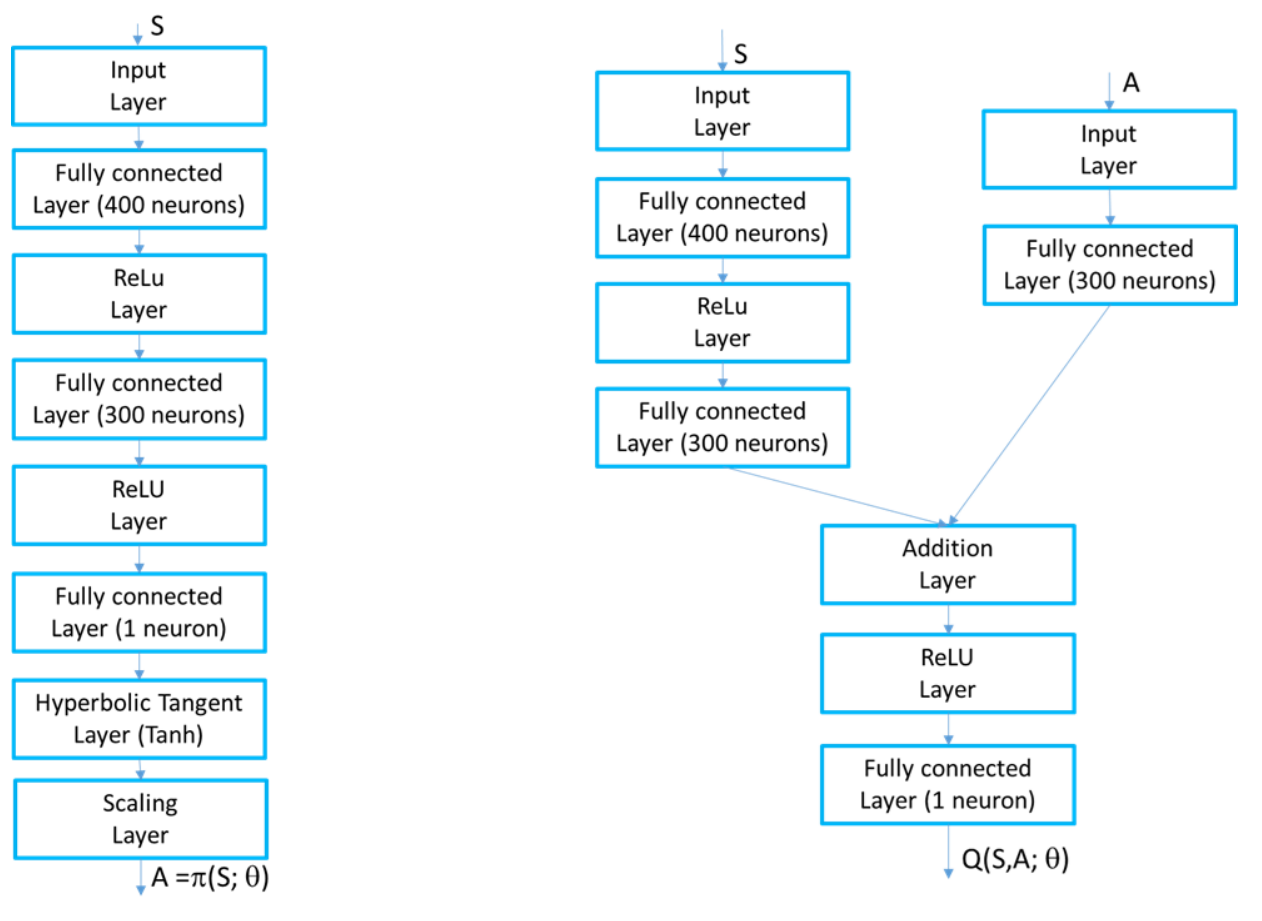}
	\caption{\label{Net} (a) Actor network (b) Critic network}
\end{figure}

When proposing an action in reinforcement learning, the dilemma of exploration and exploitation poses a great challenge. Exploration tries new, random or uncertain actions to discover new strategies while exploitation believes in the current network and utilizes it to generate the action that is assumed to achieve the best performance. DDPG achieves the balance of exploration and exploitation by adding a noise with the standard deviation $\sigma$ to the actor network. It is important to select this parameter in the implementation. Fig.~\ref{Net}(a) shows the actor network which takes in the state to generate the action. Then the action out of the actor network is added by a stochastic Gaussian or Ornstein Uhlenbeck (OU) noise $N$:
\begin{equation}
A=\pi(S;\theta_A)+N
\label{action}
\end{equation}
The standard deviation of the noise reflects if the network is in either an exploration stage or an exploitation stage. A large deviation allows the network to try out or explore new options in the action space while a small deviation implies that the network is already good enough to propose an action, therefore adding little change to the good action. A typical strategy to balance the trade-off of exploration and exploitation is that in early training, the agent is given large noise to favor exploration and build a solid knowledge base. Then in later training, the agent may shift to favor exploitation by lowering the noise as it becomes more confident in the current actor network to generate the best action.

Fig.~\ref{Net}(b) shows the critic network, which uses the state and action as input to estimate the Q-value (the expected long-term reward). Unlike Q-learning, DDPG uses target networks for stability. A target value, $y$, is calculated using the reward $R$ and the output of the target actor ($\pi'$) and target critic ($Q'$) networks:
\begin{equation}
y = R + \gamma Q'(S', \pi'(S'; \theta'_{\pi}); \theta'_{Q})
\label{ddpg_target}
\end{equation}
where $\gamma \in [0,1]$ is the discount factor. The critic network is then updated by minimizing the mean squared error loss, $L$, between this target value $y$ and the critic's current output over a mini-batch of $N$ samples:
\begin{equation}
L = \frac{1}{N} \sum_i \left( y_i - Q(S_i, A_i; \theta_C) \right)^2
\label{critic_loss}
\end{equation}

\subsection{Modeling the environment: reset and step functions}
For our path finding problem, the state of the agent at time $t$ is designed as 3D including 2D position and 1D heading angle in radians:  $S_{t} = [X_t, Y_t, \phi_t]^\top$, where $X_t$ and $Y_t$ are the coordinates of the agent at time $t$, and $\phi_t$ is the heading of the agent at time $t$. The action $A_t$ is the change in heading $\Delta \phi_t$ at time $t$. In the Matlab implementation of the DDPG learning, the reset function defines how the agent starts when a training episode ends, while the step function specifies how the state of the agent evolves for an action, how the reward is collected, and how the environment reacts to the action. In the reset function, the agent starts with the initial position and heading, but keeps the recent network parameters. Assuming a constant velocity with the prior known speed $|V|=200$, in the step function, the update of the state is given by:
\begin{equation}
\phi_{t+1} = \phi_t + A_t,
\label{phi}
\end{equation}
and
\begin{equation}
[X_{t+1},Y_{t+1}]^\top = [X_{t},Y_{t}]^\top + T_s |V| [\cos(\phi_{t+1}),\sin(\phi_{t+1})]^\top,
\label{xUpdated}
\end{equation}
where $T_s$ is the time step, set as $T_s=0.25$ for our case.

In addition, the conditions for ending the current training episode are also defined in the step function. For example, an episode can stop when the maximum number of steps within the episode is reached, the statistically evaluated reward reaches a threshold, the agent enters the no-go zone, or the agent is out of boundary.

\subsection{Reward function}
A critical part of using the DDPG machine learning method is the proper setup of the reward function. In this preliminary research, three components are used in the reward function. First, an attractive artificial potential field (PF) is used to guide the agent to move in the direction to the target final destination $P_f$. Unlike the traditional optimal control method where the target final destination is employed as the equality condition in the optimization, reaching $P_f$ in reinforcement learning is given the maximum reward. This idea was first proposed in 1984 by Hogan \cite{hogan1984} for robot manipulation. The idea has been used in the reward function for reinforcement learning in \cite{yao2020path} and for model predictive control in \cite{zmhzy23,pjzwys23}. Clearly, by taking a negative sign, we can introduce a repulsive artificial potential field to prevent the agent from entering the no-go zone. Unlike the traditional optimal control method where the no-go zone is the nonlinear path constraint in the optimization, the continuous repulsive field repels the agent from closely approaching the no-go zone. Given the desired final destination $P_f=(-200,-400)$ and the center of the no-go zone $O=(0,0)$, denote $P$ as the current agent position, we use a normalized bivariate Gaussian distribution $N(P,P_f,\sigma)$ for the attractive potential field (it rewards the agent for moving toward the destination) and a normalized bivariate distribution $N(P,O,\sigma)$ for the repulsive potential field (it penalizes the agent for moving toward the no-go zone), where $\sigma$ is the standard deviation. To reduce energy consumption, a penalty in the reward function is introduced for heading changes $A$ at the current position (this implies that the path length is minimized). Therefore, our attractive reward function is defined as:
\begin{equation}
R(P,A) = w_1 N(P,P_f,\sigma) - w_2 N(P,O,\sigma) - w_3 |A|^2
\label{e:reward}
\end{equation}
Note that $N(P,\mu,\sigma)$ is the normalized 2D probabilistic density function (pdf) at position $P$ with mean $\mu$ ($\mu=P_f$ or $\mu=O$) and covariance matrix $\sigma$.

\begin{figure}[htbp]
\centering
\includegraphics[width=0.44\textwidth]{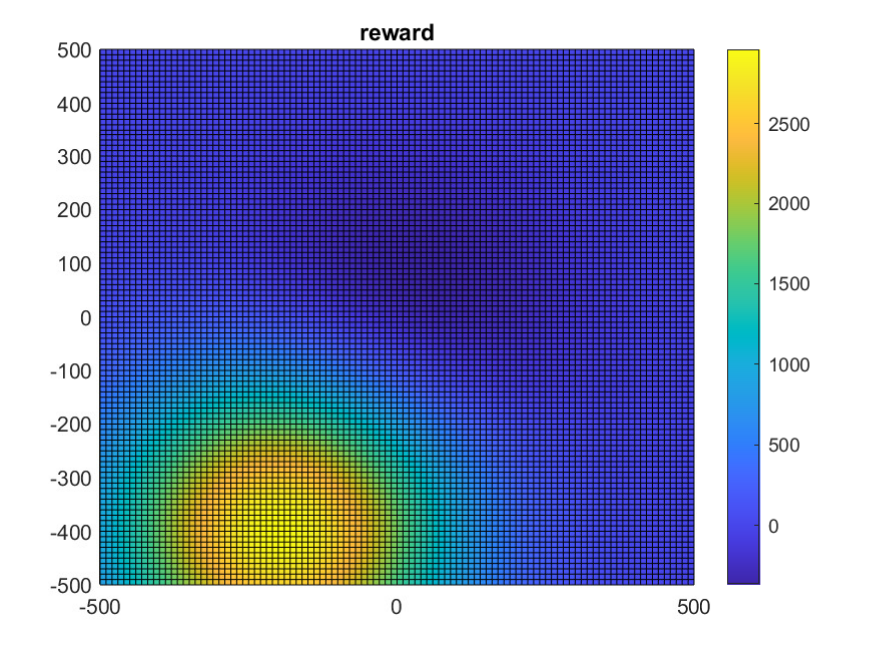}
\caption{\label{f:p1} Heat map of the reward function where $w_3$=0}
\end{figure}

The positive weights $w_1,w_2,w_3$ are carefully designed so that going around the restricted `no-go' zone has reward advantages over crossing it. In addition, the reward for approaching the target location is expected to dominate at the end of the path, and the penalty for a large control $|A|^2$ also plays a role in finding a short path. When $w_2$ is small, the agent favors the straight path regardless of the existence of the no-go zone. When $w_2$ is large, the agent has a hard time finding the target final destination. To assist the selection of the weights of $w_1$, $w_2$, $w_3$, a heat map as described in Fig.~\ref{f:p1} is used in the process. Finally, the instantaneous reward is accumulated for the long-term reward, which is optimized by the critic network according to equations (\ref{ddpg_target}) and (\ref{critic_loss}).

All hyperparameters are provided in Table~\ref{parameterTable}.
\begin{table*}[htbp]
\centering
\begin{tabular}{|c|c|c|c|}
\hline
$w_1$ & $w_2$ & $w_3$ & $\sigma$ \\
\hline
 3000 & 500 & 1.5 &
 	 $\begin{pmatrix}
 		40,000 & 100  \\
 		100 & 40,000  \\
 	\end{pmatrix}$
\\
\hline
\end{tabular}
\caption{Hyperparameters in reward function}
\label{parameterTable}
\end{table*}

\subsection{Training and testing}
\begin{figure}[htbp]
	\centering
	\includegraphics[width=0.48\textwidth]{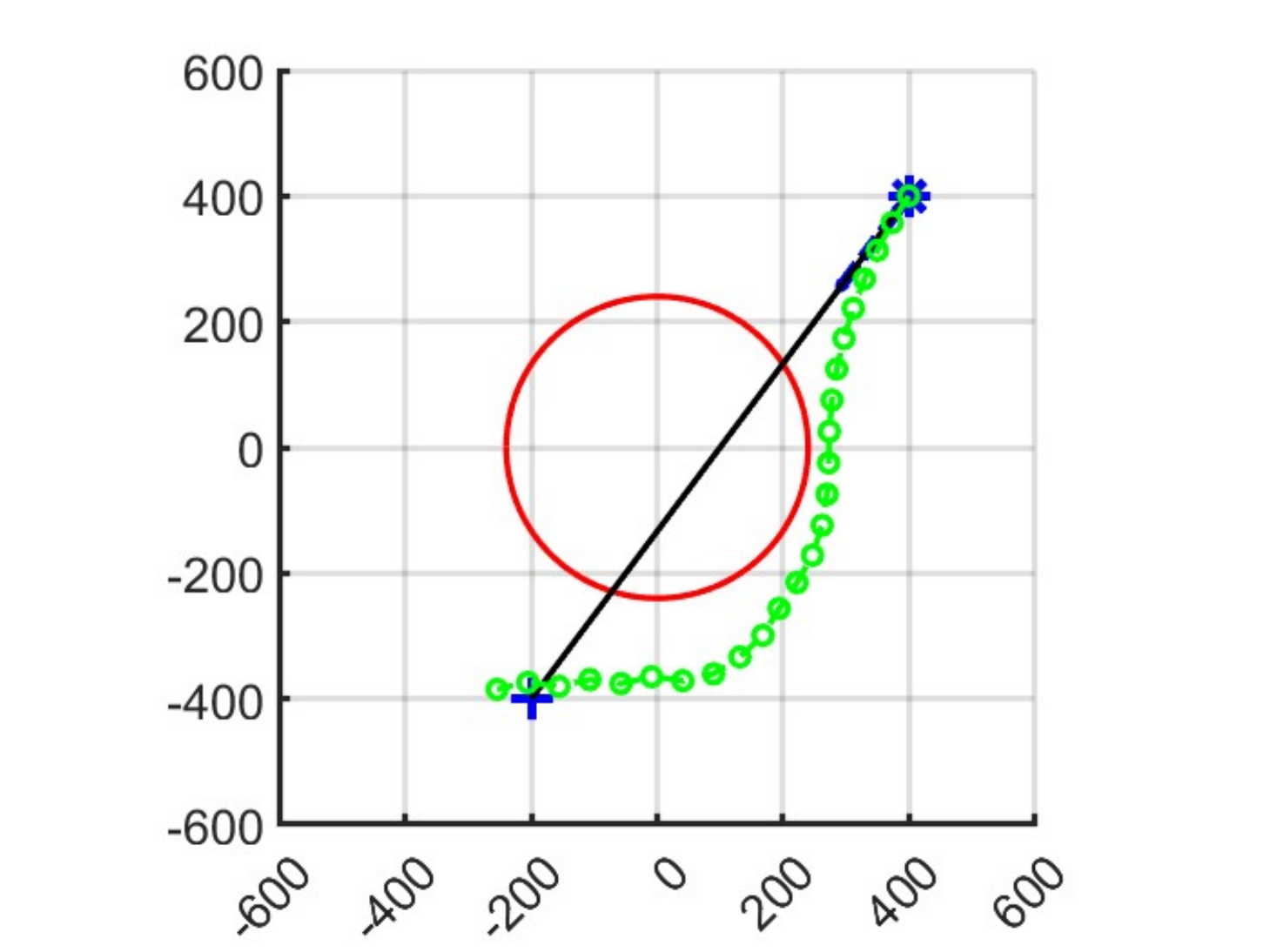}
	\includegraphics[width=0.48\textwidth]{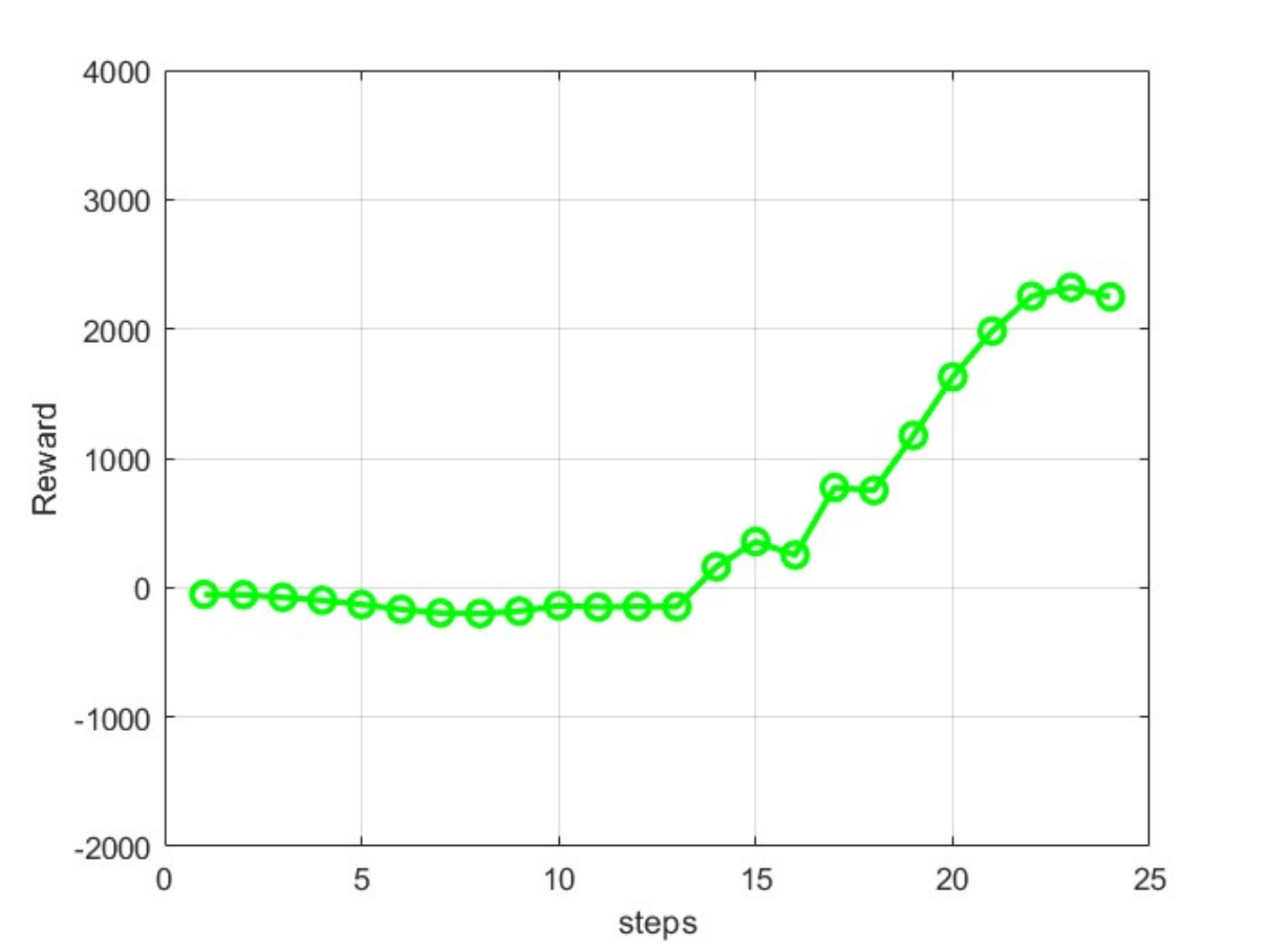}
	\includegraphics[width=0.48\textwidth]{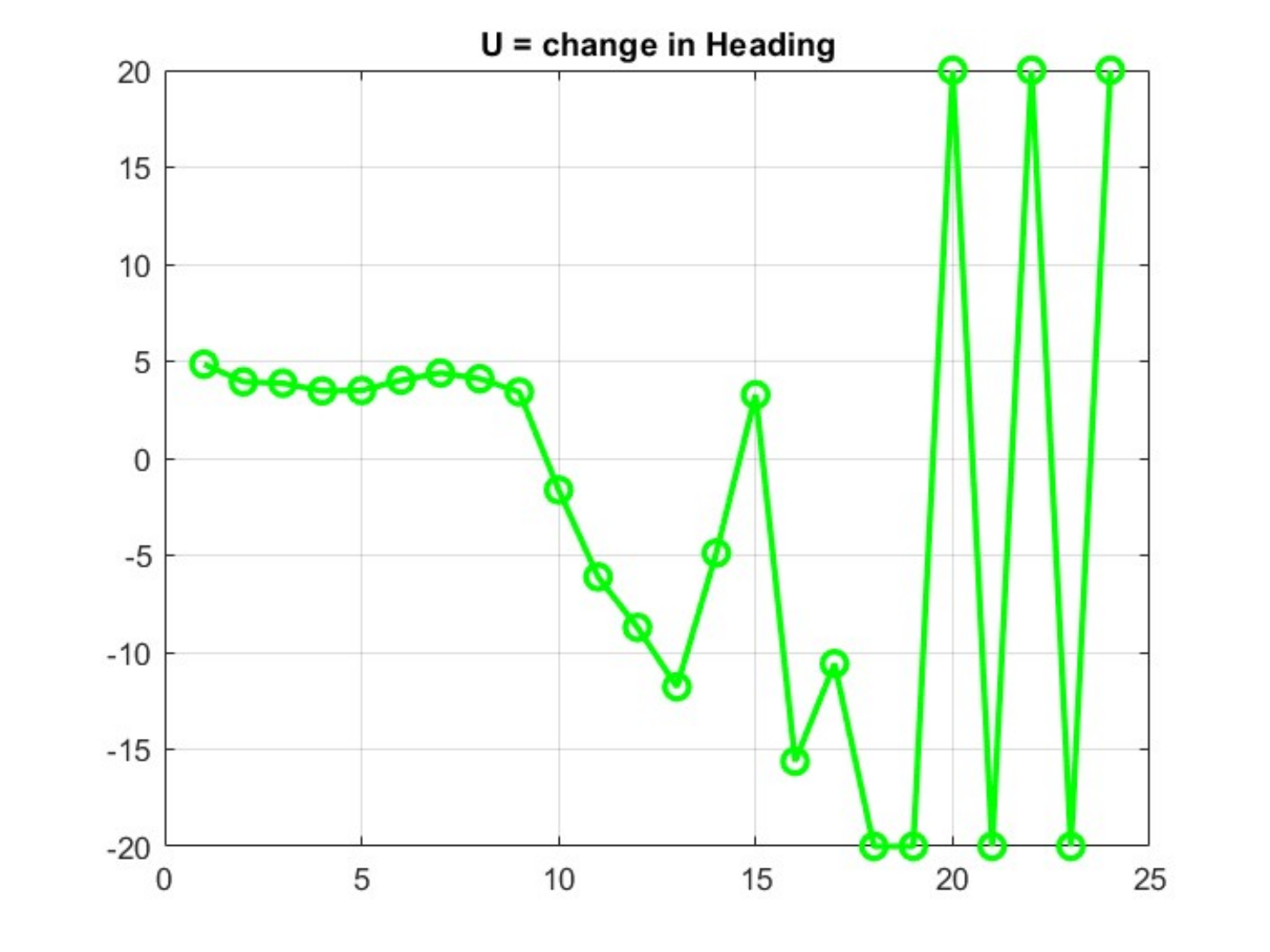}
	\caption{\label{f:training1}(a) Obtained path by DDPG with the standard deviation 5 degrees in the action noise, the resulting travel time $23 \times 0.25 = 5.75$ (b) Instant Reward (c) Action over time steps}
\end{figure}

\begin{figure}[htbp]
	\centering
	\includegraphics[width=0.48\textwidth]{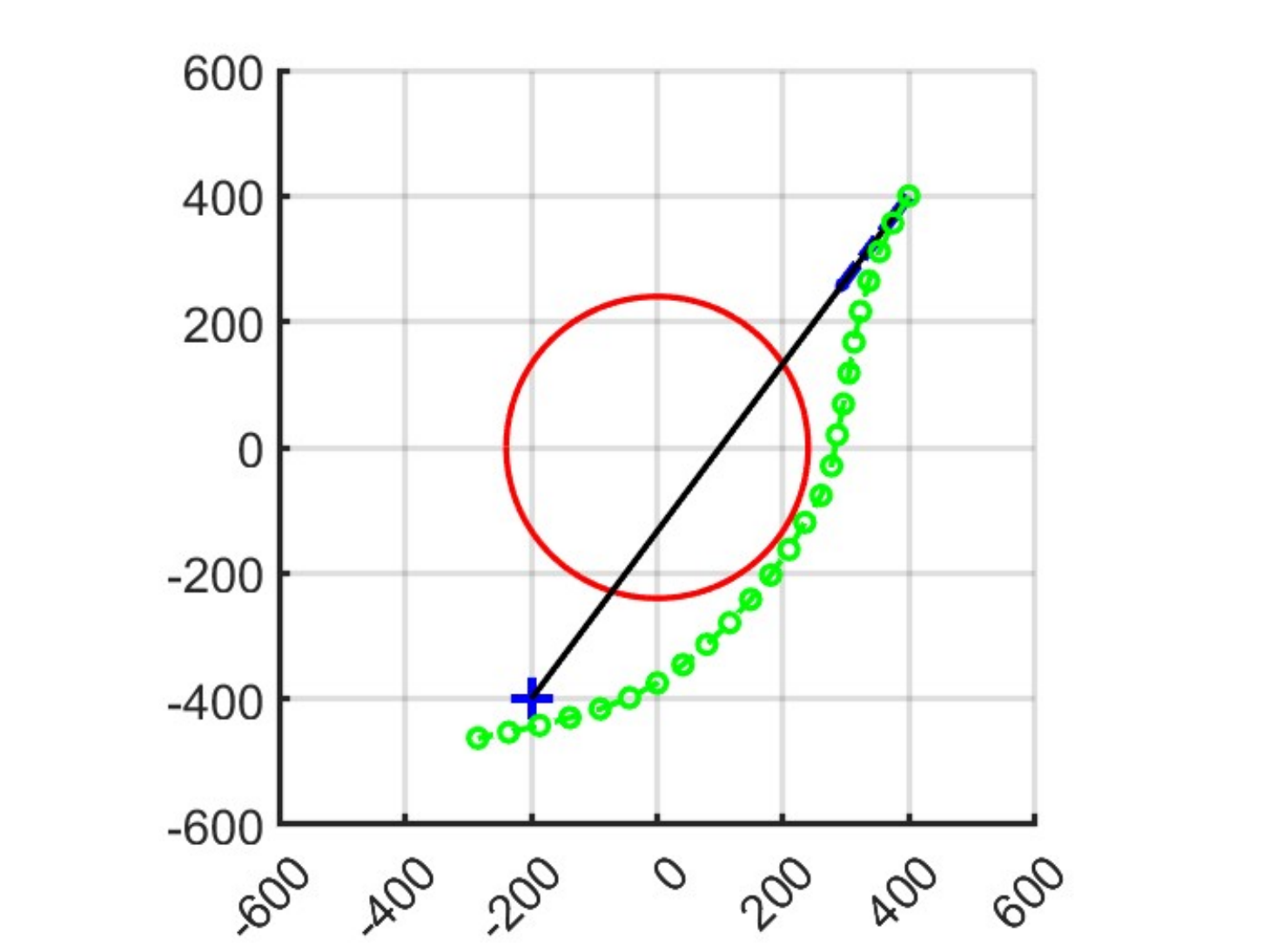}
	\includegraphics[width=0.48\textwidth]{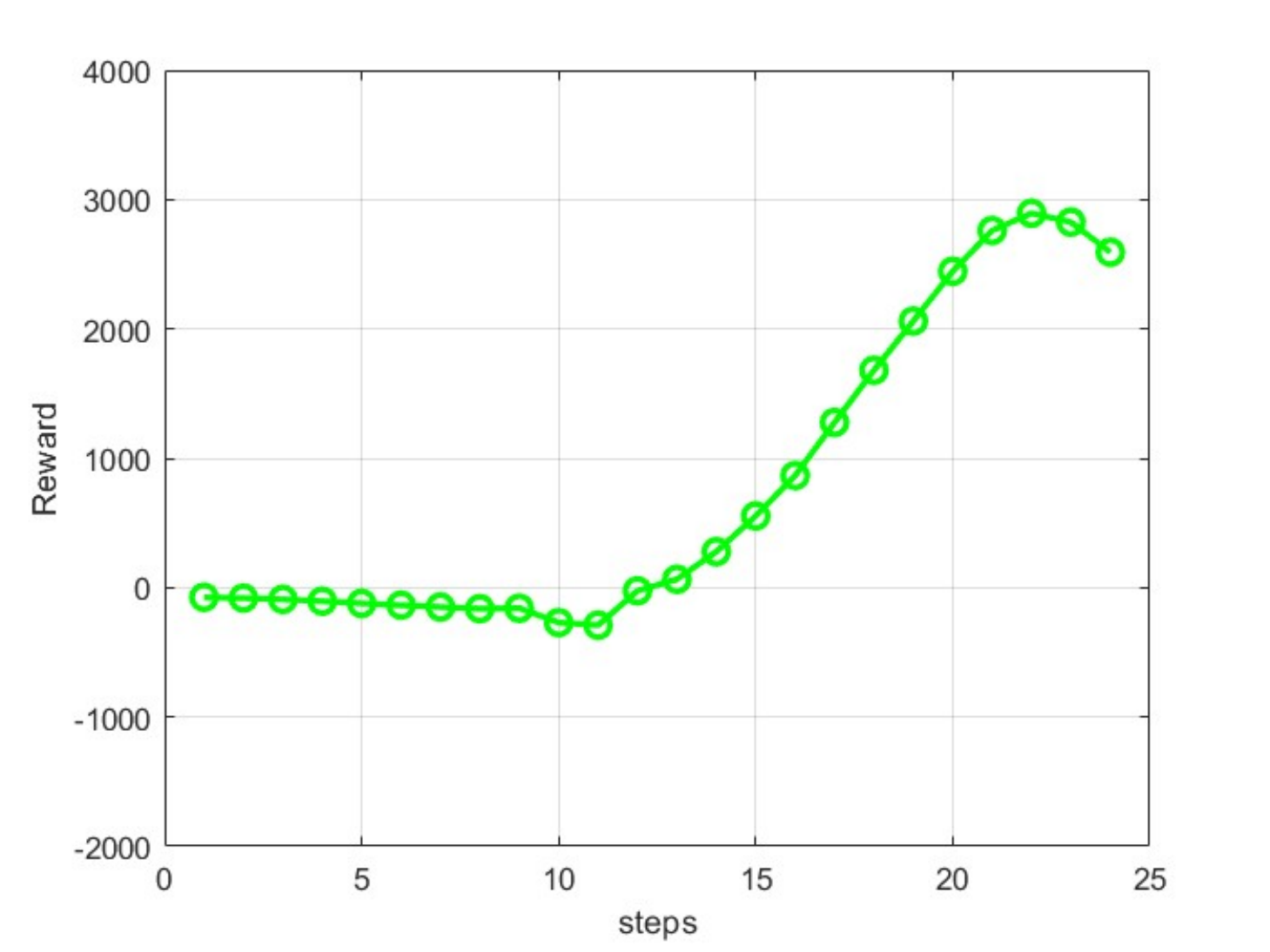}
	\includegraphics[width=0.48\textwidth]{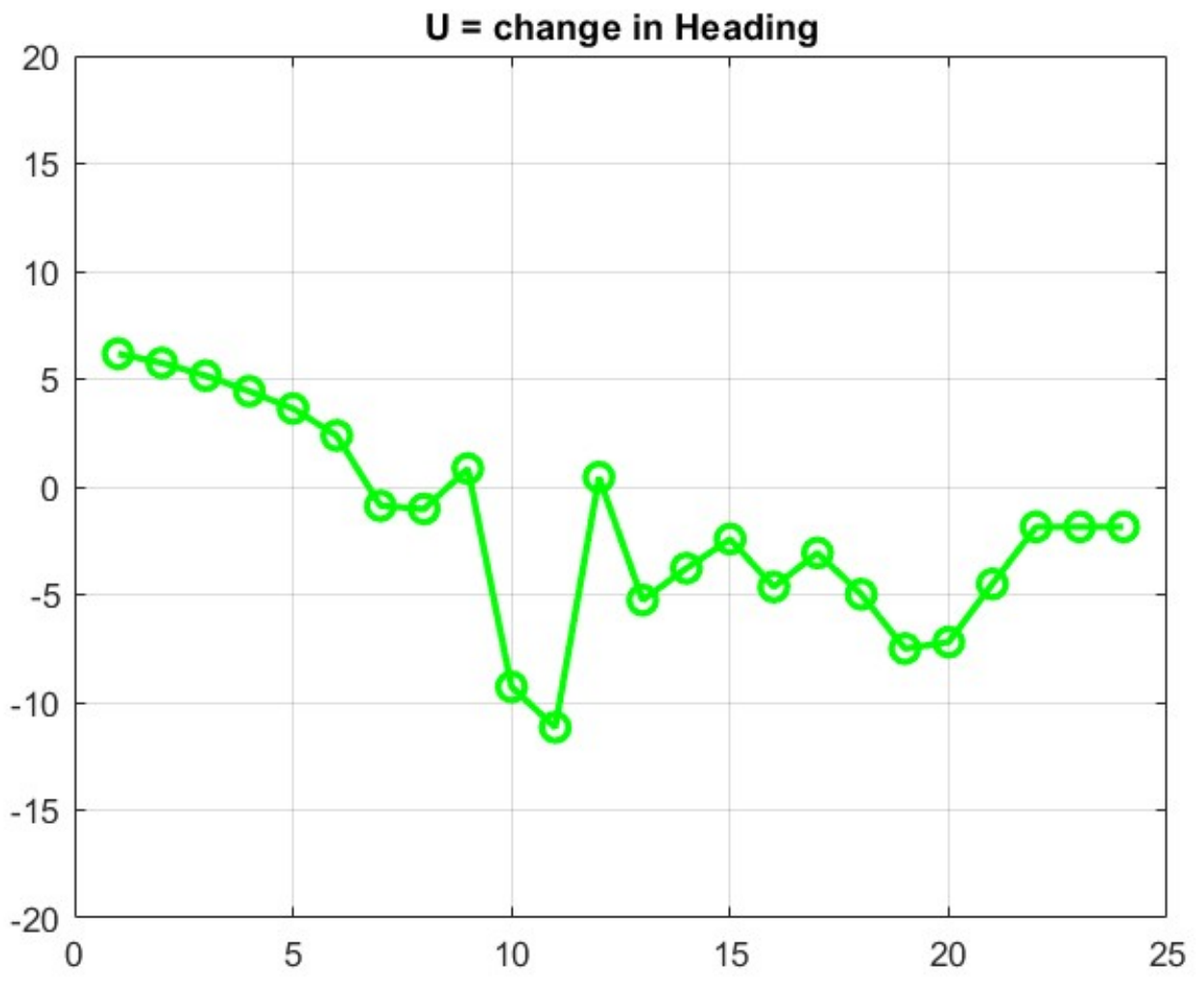}
	\caption{\label{f:training2}(a) Obtained path by further trained DDPG agent with low standard deviation in the action noise, the resulting travel time $22 \times 0.25 = 5.5$ (b) Instant Reward (c) Action over time steps}
\end{figure}

Given an initial start position $(400,400)$, the desired destination $(-200,-400)$, and the restricted `no-go' zone whose center is located at $(0,0)$ with a radius of $240$, we trained the DDPG agent and stopped the training when the statistically evaluated reward reached 12,000, and then the trained agent's performance was evaluated. The result shows that the trained agent does find a feasible path from the initial location to the desired destination. This path along with the corresponding reward and heading in each time step are plotted in Fig.~\ref{f:training1}. Direct observation of the feasible path in Fig.~\ref{f:training1} indicates that it is possible to find a shorter (better) path for this problem. Then we resumed training the agent with a smaller $\sigma$ until the accumulated reward reached 18,000. As a result, a better path is found as described in Fig.~\ref{f:training2}. Note that the maximum number of steps per episode is 24, and the agent reaches the desired destination at the 22nd step. The final travel time is computed as 5.5 seconds to reach the destination.

\begin{figure}[htbp]
\centering
\includegraphics[width=0.44\textwidth]{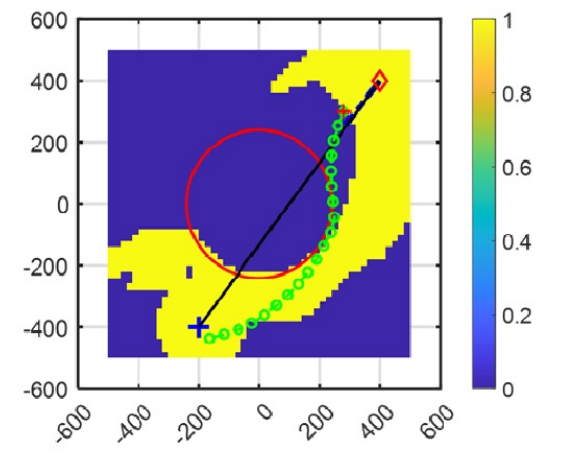}
\caption{\label{f:feasibleArea1} Feasible area of the specialized agent who is trained using $[400,400]^\top$ as the starting point}
\end{figure}

\begin{figure}[htbp]
\centering
\includegraphics[width=0.44\textwidth]{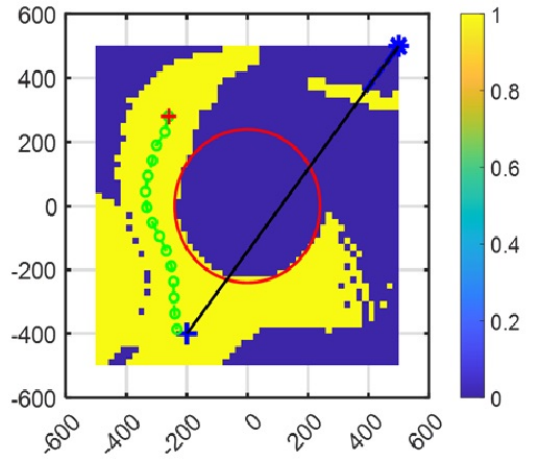}
\caption{\label{f:feasibleArea2} Feasible area of the general agent who is trained using random initial points}
\end{figure}

\begin{figure}[htbp]
	\centering
	\includegraphics[width=0.44\textwidth]{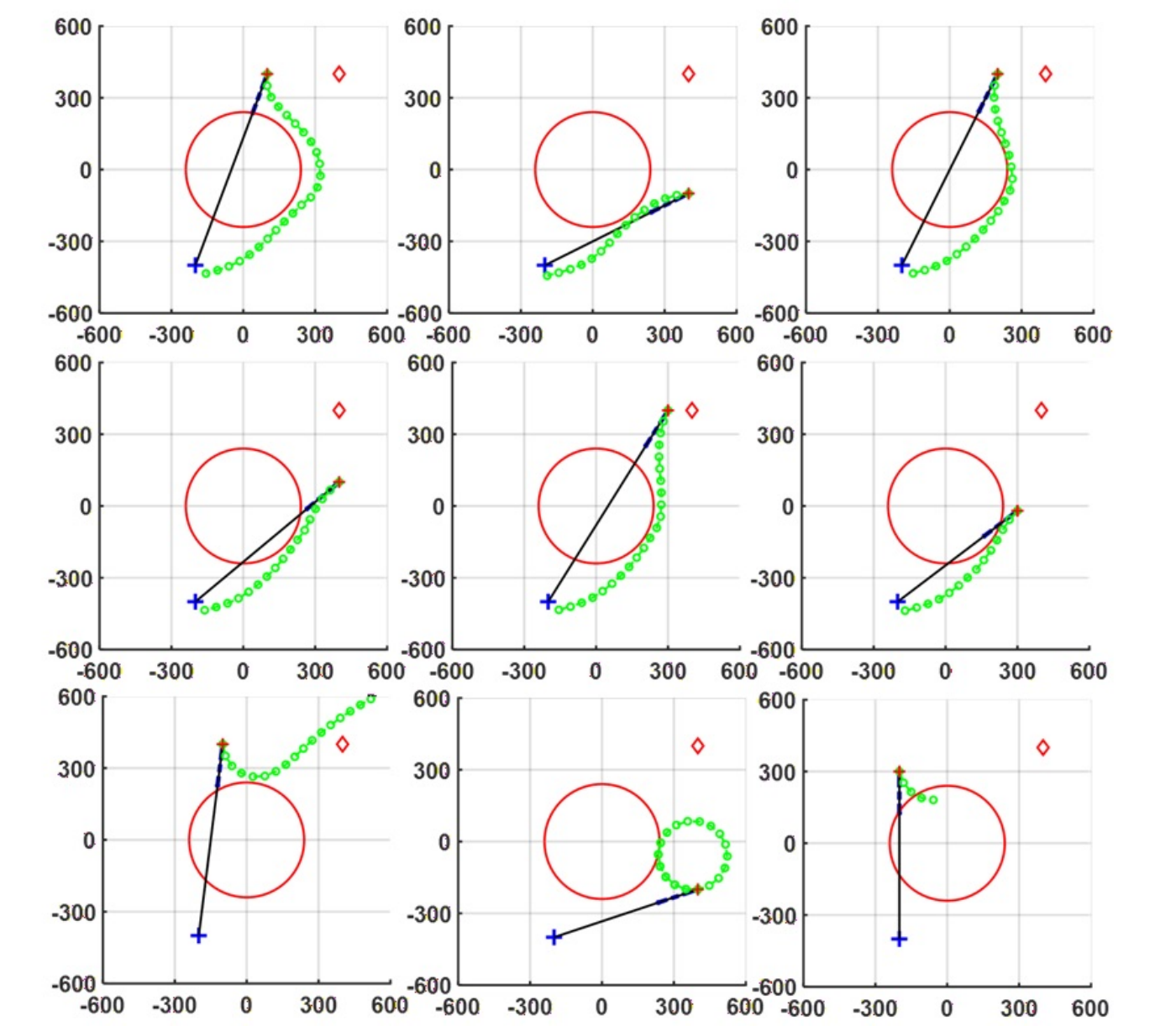}
	\caption{\label{f:ddpg1} Testing examples of the specialized DDPG agent}
\end{figure}

\begin{figure}[htbp]
	\centering
	\includegraphics[width=0.44\textwidth]{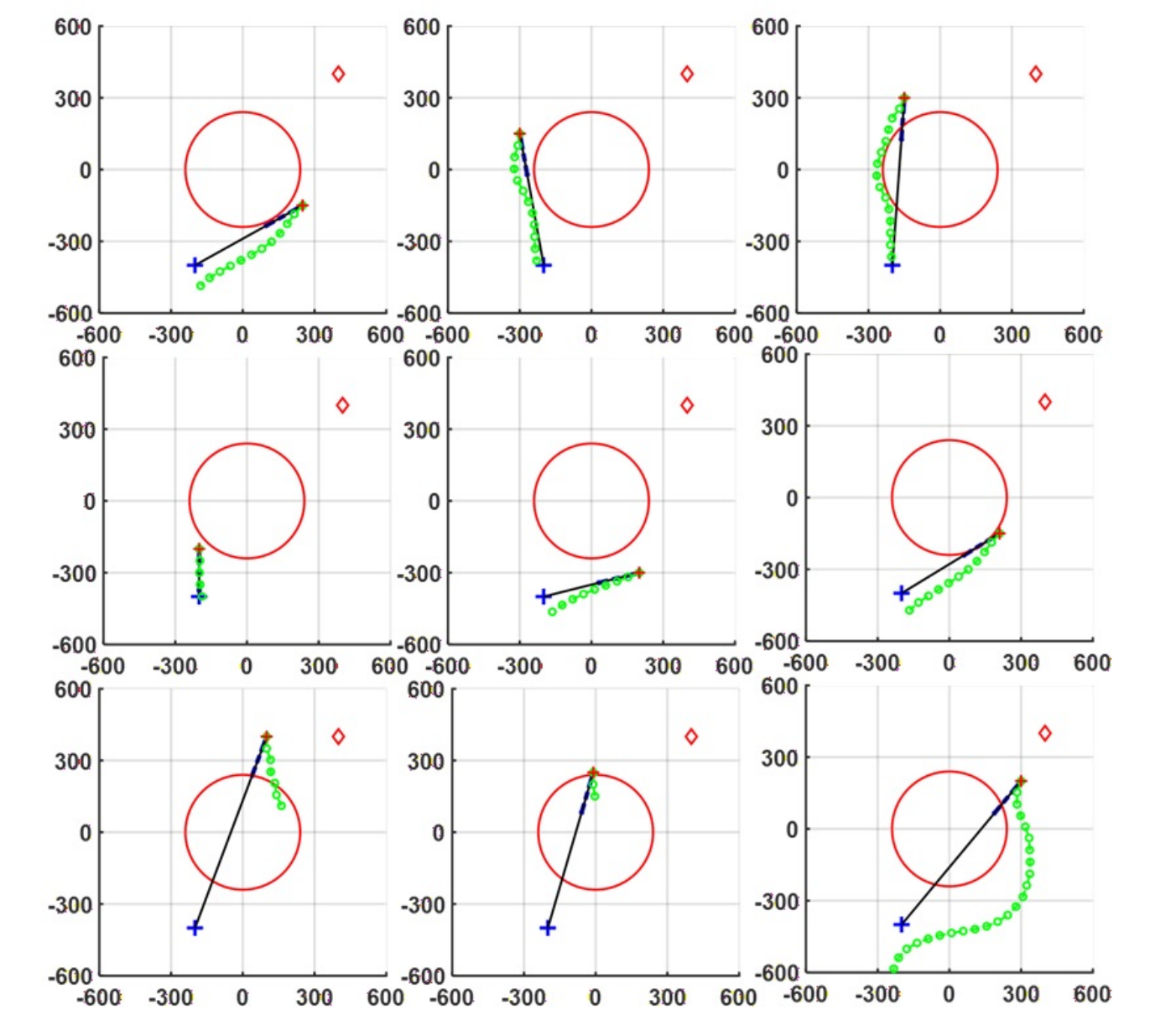}
	\caption{\label{f:ddpg2} Testing examples of the general DDPG agent}
\end{figure}

To better evaluate the performances of the agents, a feasible area (a set of starting positions) needs to be built. From any point in the feasible area, the trained agent shall provide a path to reach the proximity of the desired destination without entering the no-go zone. We train two agents to compare their feasible performance areas: a specialized agent and a general agent. The specialized agent, as mentioned above, always starts from a fixed initial point $[400,400]^\top$, while the general agent starts from randomly selected initial points during training. Once the training of the agents is done, a 1000-by-1000 grid with a grid size of 20 is used for testing. If a grid point which is the starting position can lead the agent to the destination without the agent entering the no-go zone, this grid point is marked as feasible. Fig.~\ref{f:feasibleArea1} shows the feasible area in yellow for the specialized agent, and Fig.~\ref{f:feasibleArea2} shows the feasible area for the general agent. The feasible area of the specialized agent is the neighboring area of the optimal path from the special initial point to the destination. This observation justifies the ``principle of optimality'' for using the Bellman equation in DDPG, which states that an optimal policy contains optimal sub-policies. The examples of resulting paths when testing the specialized agent are listed in Fig.~\ref{f:ddpg1}. Good paths are obtained when the starting positions are within the proximity of the specialized trained position. For the general agent, the proximity area of the destination is more likely to be feasible than the further area. In addition, the resulting paths from the general agent, as shown in Fig.~\ref{f:ddpg2}, are regarded as fair compared to those provided by the specialized agent. Future works will include the study of the reward function to improve the feasible area.

\section{Pseudo-spectral Method}\label{ps}
The objective function for a general control problem is composed of a running cost noted by the integral of the Lagrange term, and a terminating cost noted by the Mayer term:
\begin{equation}
J(S,U,t_f)=\int_{t_0}^{t_f} \textit{L}(S(t),U(t),t) \,dt +\textit{M}(t_f,S(t_f),U(t_f)),
\label{obj}
\end{equation}
where $S$ is the state vector, and $U$ is the control vector. The optimal control problem is to solve
\begin{equation}
(\hat{S},\hat{U},\hat{t}_f) = \operatorname{argmin}_{S,U,t_f} J
\label{objF}
\end{equation}
subject to the dynamic equality, path constraints, and boundary constraints, respectively:
\begin{gather}
    \frac{dS}{dt}=f_d(S(t),U(t),t), \label{diff} \\
    C(S(t),U(t)) \leq 0, \label{pathC} \\
    \textit{B}(S(t_0),S(t_f),U(t_0),U(t_f)) \leq 0. \label{boundC}
\end{gather}

Pseudo-spectral is one of the most popular numerical methods to solve the optimal control problem. Assume that LGL (Legendre-Gauss-Lobatto) collocation is applied, the multiple time interval approximation is given by:
\begin{equation}
 \min_{S,U,t_f}\sum_{i=1}^{N_I} \Bigg[ \frac{h_i}{2}\sum_{k=1}^{N_{i}}w_{i,k} \textit{L}(S_{i,k}, U_{i,k}) \Bigg] +\textit{M}(t_f,S(t_f),U(t_f))
\label{disObj}
\end{equation}
subject to continuity equality:
\begin{equation}
S_{i,N_i} = S_{i+1,1},
\label{disConst}
\end{equation}
where $i$ denotes the $i$-th mesh interval, and $N_i$ denotes the last node of the $i$-th mesh interval \cite{herber15}. For the path planning problem, the objective function becomes
\begin{equation}
(\hat{S},\hat{U},\hat{t}_f) = \operatorname{argmin}_{S,U,t_f} t_f
\label{pathObj}
\end{equation}
where $S = [P_x, P_y, V_x, V_y]^\top$ and $U = [dV_x/dt, dV_y/dt]^\top$.

Using the pseudo-spectral method, one can obtain $\hat{X}=[\hat{S}_{i,k},\hat{U}_{i,k},\hat{t}_f]^\top$ at any discrete time for the $i$th time segment and $k$th LGL (Legendre-Gauss-Lobatto) node along the path.

\begin{figure}[htbp]
	\centering
	\includegraphics[width=0.44\textwidth]{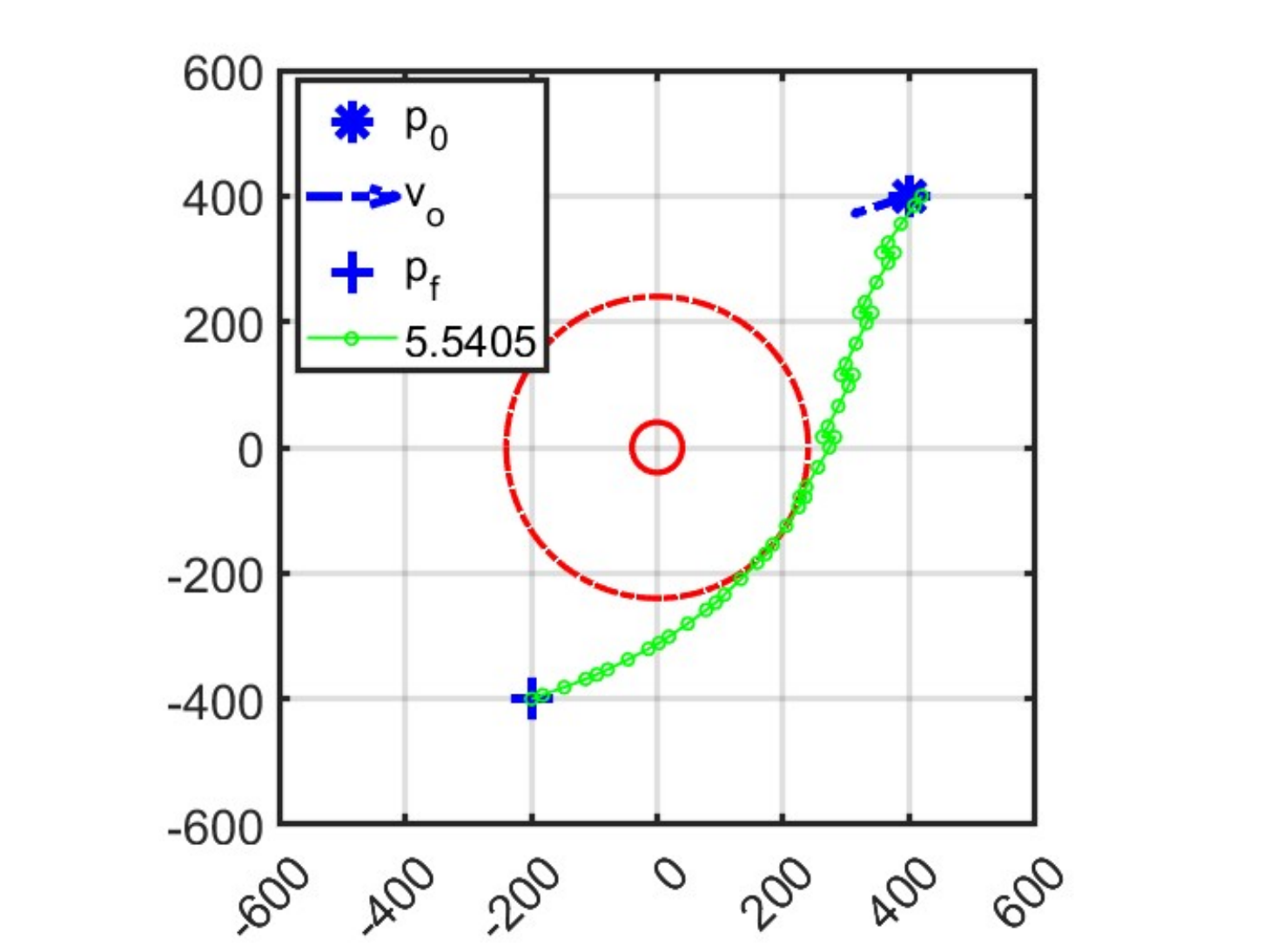}
	\caption{\label{f:training3} Obtained path by PS, Elapsed time $\approx 8$ seconds (a) constraining speed $|V|=200$, (b) no speed constraint}
\end{figure}

The pseudo-spectral method is implemented in Matlab. The code is used to solve the path planning problem with the same initial position at $(400,400)$, the same desired destination $(-200,-400)$, and the same restricted circular `no-go' zone with a center at $(0,0)$ and a radius of $240$. Constraining the agent's speed to be $200$, the resulting path from the PS method gives a 5.54 travel time in Fig.~\ref{f:training3}. The elapsed time for computing the path using the PS is around $8$ seconds on a Dell Precision 5560 (i7 Intel 2.4GHz processor, 32GB memory).

\subsection{Comparisons between reinforcement learning and optimal control}\label{Comp}
This section compares the advantages and disadvantages of reinforcement learning and optimal control (pseudo-spectral) methods when applied to the path planning problem. Setting the objective function in DDPG and PS is different for path finding. While reaching the target final destination is given the maximum reward and approaching the no-go zone is given a penalty in DDPG, the PS uses the final destination as the equality constraint and the no-go zone as nonlinear path constraints. Both DDPG and PS provide competitive paths in terms of travel time. How DDPG agents behave or perform depends on how they are trained. Although training DDPG agents takes time, the computational time of path planning when testing the DDPG agents is much less than using PS. In other words, testing a trained DDPG agent for a path solution is hundreds of times faster than using PS.

This is particularly important for real-time applications. The drawback of using the machine learning method for path planning is that it may not always provide a feasible solution, and the solution may not be optimal. While pseudo-spectral is a much slower method, if it can solve the problem, it provides an optimal solution.

\subsection{Limitations}
There are two limitations that may limit the machine learning methods from real-time applications. First, the machine learning methods may not find a feasible set that covers all starting points outside the no-go zone (we will investigate if the pseudo-spectral method can find optimal solutions for all starting points outside the no-go zone due to numerical difficulty). This drawback can be mitigated when improved reward functions are proposed. This is an ongoing effort, and we believe the goal is achievable. Second, many paths in the feasible set are not optimal; this problem may be fixed if we use the feasible paths as the initial paths and optimize the paths using a numerical optimization method. This is one of our next research topics.

\section{Conclusion and Future Work}\label{Conc}
In this paper, we developed a DDPG-based reinforcement learning method for the path planning problem to avoid a single static circular `no-go' zone, and we trained two DDPG agents: one specialized agent trained with a special initial point and one general agent trained with randomized initial points. The feasible areas are obtained for the two agents to assess where the agents can start to reach the destination. Our main purpose is to compare the performance of DDPG and that of the optimal control method (the pseudo-spectral method) to guide our future research. Our observations are (a) the machine learning method takes a longer time to train the agent, but the trained agent can provide a much faster solution than the optimal control method; (b) the pseudo-spectral method can provide an optimal solution at the expense of longer computational time, which is an obvious drawback for real-time applications. As we pointed out earlier, these are preliminary results; we are working on (a) a better training strategy to find a better reinforcement learning agent, (b) investigating if the pseudo-spectral method can find an optimal solution for all starting points outside the no-go zone, and (c) investigating if we can combine the DDPG method and the arc-search interior-point method \cite{yang25} to obtain a feasible set that covers all starting points outside of the no-go zone and to make sure that every path is optimal.

\section{Acknowledgments}
This work was supported by the Department of the Air Force under Contract No. FA955023D0001. The views and conclusions contained in this document are those of the authors and should not be interpreted as representing the official policies, either expressed or implied, of the Department of the Air Force or the U.S. Government.

\bibliography{yang}

\end{document}